\documentclass[11pt]{article}
\usepackage{mathrsfs}
\usepackage{amscd}
\usepackage{amsmath,amsfonts,amssymb,amscd}
\usepackage{indentfirst,graphics,epsfig,psfrag}
\input{epsf}
\usepackage{ifpdf}
\usepackage{enumerate}
\usepackage{appendix}
\usepackage{enumerate}
\usepackage{color}
\usepackage{lineno}
\makeatletter \@addtoreset{figure}{section} \makeatother
\makeatletter
\long\def\@makecaption#1#2{%
   \vskip 10\p@
   \setbox\@tempboxa\hbox{{#1}\ \ #2}%
   \ifdim \wd\@tempboxa >\hsize

       {#1}\ \ #2\par
   \else
       \hbox to\hsize{\hfil\box\@tempboxa\hfil}%
   \fi}
\makeatother

\newtheorem{thm}{Theorem}[section]
\newtheorem{cor}[thm]{Corollary}
\newtheorem{lem}[thm]{Lemma}

\newtheorem{pro}[thm]{Proposition}

\newcommand{\qed}{{\hfill\rule{3pt}{7pt}}}

\def\qed{\hfill \rule{4pt}{7pt}}

\begin{document}
\title{\textbf{Constructing Internally Disjoint Pendant Steiner
Trees in Cartesian Product Networks}\footnote{Supported by the
National Science Foundation of China (No. 11161037) and the Science
Found of Qinghai Province (No. 2014-ZJ-907).}}
\author{
\small  Yaping Mao\footnote{E-mail: maoyaping@ymail.com}\\[0.2cm]
\small Department of Mathematics, Qinghai Normal\\
\small University, Xining, Qinghai 810008, China\\
\small }
\date{}
\maketitle
\begin{abstract}
The concept of pedant tree-connectivity was introduced by Hager in
1985. For a graph $G=(V,E)$ and a set $S\subseteq V(G)$ of at least
two vertices, \emph{an $S$-Steiner tree} or \emph{a Steiner tree
connecting $S$} (or simply, \emph{an $S$-tree}) is a such subgraph
$T=(V',E')$ of $G$ that is a tree with $S\subseteq V'$. For an
$S$-Steiner tree, if the degree of each vertex in $S$ is equal to
one, then this tree is called a \emph{pedant $S$-Steiner tree}. Two
pedant $S$-Steiner trees $T$ and $T'$ are said to be
\emph{internally disjoint} if $E(T)\cap E(T')=\varnothing$ and
$V(T)\cap V(T')=S$. For $S\subseteq V(G)$ and $|S|\geq 2$, the
\emph{local pedant tree-connectivity} $\tau_G(S)$ is the maximum
number of internally disjoint pedant $S$-Steiner trees in $G$. For
an integer $k$ with $2\leq k\leq n$, \emph{pedant tree
$k$-connectivity} is defined as
$\tau_k(G)=\min\{\tau_G(S)\,|\,S\subseteq V(G),|S|=k\}$. In this
paper, we prove that for any two connected graphs $G$ and $H$,
$\tau_3(G\Box H)\geq \min\{3\lfloor\frac{\tau_3(G)}{2}\rfloor,3\lfloor\frac{\tau_3(H)}{2}\rfloor\}$. Moreover,
the bound is sharp.\\[2mm]
\end{abstract}

\section{Introduction}

A processor network is expressed as a graph, where a node is a
processor and an edge is a communication link. Broadcasting is the
process of sending a message from the source node to all other nodes
in a network. It can be accomplished by message dissemination in
such a way that each node repeatedly receives and forwards messages.
Some of the nodes and/or links may be faulty. However, multiple
copies of messages can be disseminated through disjoint paths. We
say that the broadcasting succeeds if all the healthy nodes in the
network finally obtain the correct message from the source node
within a certain limit of time. A lot of attention has been devoted
to fault-tolerant broadcasting in networks \cite{Fragopoulou,
Hedetniemi, Jalote, Ramanathan}. In order to measure the ability of
fault-tolerance, the above path structure connecting two nodes are
generalized into some tree structures connecting more than two nods,
see \cite{Ku, LLSun, LM1}. To show these generalizations clearly, we
must state from the connectivity in Graph Theory.

All graphs considered in this paper are undirected, finite and
simple. We refer to the book \cite{bondy} for graph theoretical
notation and terminology not described here. For a graph $G$, let
$V(G)$, $E(G)$ and $\delta(G)$ denote the set of vertices, the set
of edges and the minimum degree of $G$, respectively. Connectivity
is one of the most basic concepts of graph-theoretic subjects, both
in combinatorial sense and the algorithmic sense. It is well-known
that the classical connectivity has two equivalent definitions. The
\emph{connectivity} of $G$, written $\kappa(G)$, is the minimum
order of a vertex set $S\subseteq V(G)$ such that $G\setminus S$ is
disconnected or has only one vertex. We call this definition the
`cut' version definition of connectivity. A well-known theorem of
Whitney provides an equivalent definition of connectivity, which can
be called the `path' version definition of connectivity. For any two
distinct vertices $x$ and $y$ in $G$, the \emph{local connectivity}
$\kappa_{G}(x,y)$ is the maximum number of internally disjoint paths
connecting $x$ and $y$. Then
$\kappa(G)=\min\{\kappa_{G}(x,y)\,|\,x,y\in V(G),x\neq y\}$ is
defined to be the \emph{connectivity} of $G$. For connectivity,
Oellermann gave a survey paper on this subject; see
\cite{Oellermann2}.

Although there are many elegant and powerful results on connectivity
in Graph Theory, the basic notation of classical connectivity may
not be general enough to capture some computational settings. So
people want to generalize this concept. For the `cut' version
definition of connectivity, we find the above minimum vertex set
without regard the number of components of $G\setminus S$. Two
graphs with the same connectivity may have differing degrees of
vulnerability in the sense that the deletion of a vertex cut-set of
minimum cardinality from one graph may produce a graph with
considerably more components than in the case of the other graph.
For example, the star $K_{1,n}$ and the path $P_{n+1}\ (n\geq 3)$
are both trees of order $n+1$ and therefore connectivity $1$, but
the deletion of a cut-vertex from $K_{1,n}$ produces a graph with
$n$ components while the deletion of a cut-vertex from $P_{n+1}$
produces only two components. Chartrand et al. \cite{Chartrand1}
generalized the `cut' version definition of connectivity. For an
integer $k \ (k\geq 2)$ and a graph $G$ of order $n \ (n\geq k)$,
the \emph{$k$-connectivity} $\kappa'_k(G)$ is the smallest number of
vertices whose removal from $G$ of order $n \ (n\geq k)$ produces a
graph with at least $k$ components or a graph with fewer than $k$
vertices. Thus, for $k=2$, $\kappa'_2(G)=\kappa(G)$. For more
details about $k$-connectivity, we refer to \cite{Chartrand1, Day,
Oellermann2, Oellermann3}.

The generalized connectivity of a graph $G$, introduced by Hager \cite{Hager}, is
a natural generalization of the `path' version definition of
connectivity. For a graph $G=(V,E)$ and a set $S\subseteq V(G)$ of
at least two vertices, \emph{an $S$-Steiner tree} or \emph{a Steiner
tree connecting $S$} (or simply, \emph{an $S$-tree}) is a such
subgraph $T=(V',E')$ of $G$ that is a tree with $S\subseteq V'$.
Note that when $|S|=2$ an $S$-Steiner tree is just a path
connecting the two vertices of $S$. Two $S$-Steiner trees $T$ and $T'$
are said to be \emph{internally disjoint} if
$E(T)\cap E(T')=\varnothing$ and $V(T)\cap V(T')=S$. For $S\subseteq
V(G)$ and $|S|\geq 2$, the \emph{generalized local connectivity}
$\kappa_G(S)$ is the maximum number of internally disjoint $S$-Steiner trees
in $G$, that is, we search for the maximum
cardinality of edge-disjoint trees which include $S$ and are vertex
disjoint with the exception of $S$. For an integer $k$ with $2\leq
k\leq n$, \emph{generalized $k$-connectivity} (or
\emph{$k$-tree-connectivity}) is defined as
$\kappa_k(G)=\min\{\kappa_G(S)\,|\,S\subseteq V(G),|S|=k\}$, that
is, $\kappa_k(G)$ is the minimum value of $\kappa_G(S)$ when $S$
runs over all $k$-subsets of $V(G)$. Clearly, when $|S|=2$,
$\kappa_2(G)$ is nothing new but the connectivity $\kappa(G)$ of
$G$, that is, $\kappa_2(G)=\kappa(G)$, which is the reason why one
addresses $\kappa_k(G)$ as the generalized connectivity of $G$. By
convention, for a connected graph $G$ with less than $k$ vertices,
we set $\kappa_k(G)=1$. Set $\kappa_k(G)=0$ when $G$ is
disconnected. Note that the generalized $k$-connectivity and
$k$-connectivity of a graph are indeed different. Take for example,
the graph $H_1$ obtained from a triangle with vertex set
$\{v_1,v_2,v_3\}$ by adding three new vertices $u_1,u_2,u_3$ and
joining $v_i$ to $u_i$ by an edge for $1 \leq i\leq 3$. Then
$\kappa_3(H_1)=1$ but $\kappa'_3(H_1)=2$. There are many results on
the generalized connectivity, see \cite{Chartrand2,
LLSun, LL, LLZ, LM1, LM2, LM3, LM4, LMS, Okamoto}.

The concept of pedant-tree connectivity \cite{Hager} was introduced
by Hager in 1985, which is specialization of generalized
connectivity (or \emph{$k$-tree-connectivity}) but a generalization
of classical connectivity. For an $S$-Steiner tree, if the degree of
each vertex in $S$ is equal to one, then this tree is called a
\emph{pedant $S$-Steiner tree}. Two pedant $S$-Steiner trees $T$ and
$T'$ are said to be \emph{internally disjoint} if $E(T)\cap
E(T')=\varnothing$ and $V(T)\cap V(T')=S$. For $S\subseteq V(G)$ and
$|S|\geq 2$, the \emph{local pedant-tree connectivity} $\tau_G(S)$
is the maximum number of internally disjoint pedant $S$-Steiner
trees in $G$. For an integer $k$ with $2\leq k\leq n$,
\emph{pedant-tree $k$-connectivity} is defined as
$\tau_k(G)=\min\{\tau_G(S)\,|\,S\subseteq V(G),|S|=k\}$. Set
$\kappa_k(G)=0$ when $G$ is disconnected. When $k=2$,
$\tau_2(G)=\tau(G)$ is just the connectivity of a graph $G$.

In {\upshape\cite{Hager}}, Hager derived the following results.

\begin{lem}{\upshape\cite{Hager}}\label{lem1-1}
Let $G$ be a graph. If $\tau_k(G)\geq \ell$, then $\delta(G)\geq
k+\ell-1$.
\end{lem}

\begin{lem}{\upshape\cite{Hager}}\label{lem1-2}
Let $G$ be a graph. If $\tau_k(G)\geq \ell$, then $\kappa(G)\geq
k+\ell-2$.
\end{lem}

Li et al. {\upshape\cite{LLZ}} obtained the following result.

\begin{lem}{\upshape\cite{LLZ}}\label{lem1-3}
Let $G$ be a connected graph with minimum degree $\delta$. If there are two adjacent vertices of degree $\delta$, then
$\kappa(G)\leq \delta(G)-1$.
\end{lem}

It is clear that generalized $k$-connectivity (or
\emph{$k$-tree-connectivity}) and pedant-tree $k$-connectivity of a
graph are indeed different. For example, let $H=W_n$ be a wheel of
order $n$. From Lemma \ref{lem1-1}, we have $\tau_3(H)\leq 1$. One
can check that for any $S\subseteq V(H)$ with $|S|=3$,
$\tau_{H}(S)\geq 1$. Therefore, $\tau_3(H)=1$. From Lemma
\ref{lem1-3}, we have $\kappa_3(H)\leq \delta(H)-1=3-1=2$. One can
check that for any $S\subseteq V(G)$ with $|S|=3$, $\kappa_{H}(S)\geq
2$. Therefore, $\kappa_3(H)=2$.

In fact, Mader \cite{Mader3} studied an extension of Menger's
theorem to independent sets of three or more vertices. We know that
from Menger's theorem that if $S=\{u,v\}$ is a set of two
independent vertices in a graph $G$, then the maximum number of
internally disjoint $u$-$v$ paths in $G$ equals the minimum number
of vertices that separate $u$ and $v$. For a set
$S=\{u_1,u_2,\cdots,u_k\}$ of $k \ (k\geq 2)$ vertices in a graph
$G$, an \emph{$S$-path} is defined as a path between a pair of
vertices of $S$ that contains no other vertices of $S$. Two
$S$-paths $P_1$ and $P_2$ are said to be \emph{internally disjoint}
if they are vertex-disjoint except for the vertices of $S$. If $S$
is a set of independent vertices of a graph $G$, then a vertex set
$U\subseteq V(G)$ with $U\cap S=\varnothing$ is said to
\emph{totally separate $S$} if every two vertices of $S$ belong to
different components of $G\setminus U$. Let $S$ be a set of at least
three independent vertices in a graph $G$. Let $\mu(G)$ denote the
maximum number of internally disjoint $S$-paths and $\mu'(G)$ the
minimum number of vertices that totally separate $S$. A natural
extension of Menger' s theorem may well be suggested, namely: If $S$
is a set of independent vertices of a graph $G$ and $|S|\geq 3$,
then $\mu(S)=\mu'(S)$. However, the statement is not true in
general. Take the above graph $G_0$ for example. For
$S=\{v_1,v_2,v_3\}$, $\mu(S)=1$ but $\mu'(S)=2$. Mader proved that
$\mu(S)\geq \frac{1}{2}\mu'(S)$. Moreover, the bound is sharp.
Lov\'{a}sz conjectured an edge analogue of this result and Mader
proved this conjecture and established its sharpness. For more
details, we refer to \cite{Mader3, Mader4, Oellermann1}.

In addition to being a natural combinatorial measure, pedant tree $k$-connectivity and
generalized $k$-connectivity can be motivated by its interesting
interpretation in practice. For example, suppose that $G$ represents
a network. If one considers to connect a pair of vertices of $G$,
then a path is used to connect them. However, if one wants to
connect a set $S$ of vertices of $G$ with $|S|\geq 3$, then a tree
has to be used to connect them. This kind of tree for connecting a
set of vertices is usually called a {\it Steiner tree}, and
popularly used in the physical design of VLSI circuits (see
\cite{Grotschel1, Grotschel2, Sherwani}). In this application, a
Steiner tree is needed to share an electric signal by a set of
terminal nodes. Steiner tree is also used in computer communication
networks (see \cite{Du}) and optical wireless communication networks
(see \cite{Cheng}). Usually, one wants to consider how tough a
network can be, for the connection of a set of vertices. Then, the
number of totally independent ways to connect them is a measure for
this purpose. The generalized $k$-connectivity can serve for
measuring the capability of a network $G$ to connect any $k$
vertices in $G$.

Product networks were proposed based upon the idea of using the
cross product as a tool for ``combining'' two known graphs with
established properties to obtain a new one that inherits properties
from both \cite{DayA}. Recently, there has been an increasing
interest in a class of interconnection networks called Cartesian
product networks; see \cite{Bao, DayA, Ku, LLSun}.

The \emph{Cartesian product} of two graphs $G$ and $H$, written as
$G\Box H$, is the graph with vertex set $V(G)\times V(H)$, in which
two vertices $(u,v)$ and $(u',v')$ are adjacent if and only if
$u=u'$ and $(v,v')\in E(H)$, or $v=v'$ and $(u,u')\in E(G)$.

In this paper, we obtain the following lower bound of
$\tau_3(G\Box H)$.

\begin{thm}\label{th1-4}
Let $G$ and $H$ be two connected graphs. Then
$$
\tau_3(G\Box H)\geq \min\left\{3\left\lfloor\frac{\tau_3(G)}{2}\right\rfloor,3\left\lfloor\frac{\tau_3(H)}{2}\right\rfloor\right\}.
$$
Moreover, the bound is sharp.
\end{thm}

\section{Proof of Theorem \ref{th1-4}}

In this section, let $G$ and $H$ be two connected graphs with
$V(G)=\{u_1,u_2,\ldots,u_{n}\}$ and $V(H)=\{v_1,v_2,\ldots,v_{m}\}$,
respectively. Then $V(G\Box H)=\{(u_i,v_j)\,|\,1\leq i\leq n, \
1\leq j\leq m\}$. For $v\in V(H)$, we use $G(v)$ to denote the
subgraph of $G\Box H$ induced by the vertex set
$\{(u_i,v)\,|\,1\leq i\leq n\}$. Similarly, for $u\in V(G)$, we use
$H(u)$ to denote the subgraph of $G\Box H$ induced by the vertex
set $\{(u,v_j)\,|\,1\leq j\leq m\}$. In the sequel, let $K_{s,t}$,
$K_{n}$ and $P_n$ denote the complete bipartite graph of order
$s+t$, complete graph of order $n$, and path of order $n$,
respectively. If $G$ is a connected graph and $x,y\in V(G)$, then
the \emph{distance} $d_G(x,y)$ between $x$ and $y$ is the length of
a shortest path connecting $x$ and $y$ in $G$.

We now introduce the general idea of the proof of Theorem
\ref{th1-4}, with a running example (corresponding to Fig. 1). From
the definition, Cartesian product graph $G\Box H$ is a graph
obtained by replacing each vertex of $G$ by a copy of $H$ and
replacing each edge of $G$ by a perfect matching of a complete
bipartite graph $K_{m,m}$. Recall that
$V(G)=\{u_1,u_2,\ldots,u_{n}\}$. Clearly, $V(G\Box
H)=\bigcup_{i=1}^nV(H(u_i))$. Take for example, let $G=K_8$ (see
Fig. 1 $(a)$). Set $V(K_8)=\{u_i\,|\,1\leq i\leq 8\}$ and
$|V(H)|=m$. Then $K_8\Box H$ is a graph obtained by replacing each
vertex of $K_8$ by a copy of $H$ and replacing each edge of $K_8$ by
a perfect matching of complete bipartite graph $K_{m,m}$ (see Fig. 1
$(e)$). Clearly, $V(K_8\Box H)=\bigcup_{i=1}^8V(H(u_i))$ (see Fig. 1
$(e)$).

In this section, we give the proof of Theorem \ref{th1-4}. For two
connected graphs $G$ and $H$, we prove that $\tau_3(G\Box H)\geq
\min\{3\lfloor\frac{\tau_3(G)}{2}\rfloor,3\lfloor\frac{\tau_3(H)}{2}\rfloor\}$.
By the symmetry of Cartesian product graphs, we assume
$\tau_3(H)\geq \tau_3(G)$. We need to show $\tau_3(G\Box H)\geq
3\lfloor\frac{\tau_3(G)}{2}\rfloor$. Set $\tau_3(G)=k$ and
$\tau_3(H)=\ell$. From the definition of $\tau_3(G\Box H)$, it
suffices to show that $\kappa_{G\Box H}(S)\geq
3\lfloor\frac{k}{2}\rfloor$ for any $S\subseteq V(G\Box H)$ and
$|S|=3$. Furthermore, from the definition of $\kappa_{G\Box H}(S)$,
we need to find out $3\lfloor\frac{k}{2}\rfloor$ internally disjoint
pedant $S$-Steiner trees in $G\Box H$. Let $S=\{x,y,z\}$. Recall
that $V(G)=\{u_1,u_2,\ldots,u_{n}\}$. From the above analysis, we
know that $x,y,z\in V(G\Box H)=\bigcup_{i=1}^nV(H(u_i))$. Without
loss of generality, let $x\in H(u_i)$, $y\in H(u_j)$ and $z\in
H(u_k)$ (note that $u_i,u_j,u_k$ are not necessarily different). For
the above example, we have $x,y,z\in V(K_8\Box
H)=\bigcup_{i=1}^8V(H(u_i))$. Without loss of generality, let $x\in
H(u_1)$, $y\in H(u_2)$ and $z\in H(u_3)$ (see Fig. 1 $(e)$).

Because $u_i,u_j,u_k\in V(G)$ and $\tau_3(G)=k$, there are $k$
internally disjoint pedant Steiner trees connecting
$\{u_i,u_j,u_k\}$, say $T_1,T_2,\cdots,T_{k}$. Note that
$\bigcup_{i=1}^{\ell} T_i$ is a subgraph of $G$. Let $y',z'$ be the
vertices corresponding to $y,z$ in $H(u_i)$. Since $\tau_3(H)=\ell$,
there are $\ell$ internally disjoint pedant Steiner trees connecting
$\{x,y',z'\}$ in $H(u_i)$, say $T_1',T_2',\cdots,T_{\ell}'$. Thus
$(\bigcup_{i=1}^{k} T_i)\Box (\bigcup_{j=1}^{\ell} T_j')$ is a
subgraph of $G\Box H$. For the above example, we have
$\tau_3(G)=\tau_3(K_8)=k=5\leq \ell$. It suffices to prove that
$\tau_3(G\Box H)\geq
3\lfloor\frac{\tau_3(G)}{2}\rfloor=3\lfloor\frac{k}{2}\rfloor$.
Clearly, there are $k=5$ internally disjoint pedant Steiner trees
connecting $\{u_1,u_2,u_3\}$, say $T_1,T_2,T_3,T_4,T_5$ (see
$T_1,T_2,T_3,T_4$ in Fig. 1 $(b),(c)$). Note that $T_1\cup T_2$ or
$T_3\cup T_4$ is a subgraph of $G$ (see Fig. 1 $(b),(c)$). Then
$(\bigcup_{i=1}^{4} T_i)\Box (\bigcup_{j=1}^{\ell} T_j')$ is a
subgraph of $G\Box H$ (see Fig. 1 $(d),(h)$).

If we can prove that $\tau_{(\bigcup_{i=1}^{k} T_i)\Box
(\bigcup_{j=1}^{\ell} T_j')}(S)\geq 3\lfloor\frac{k}{2}\rfloor$ for
$S=\{x,y,z\}$, then $\tau_{G\Box H}(S)\geq \tau_{(\bigcup_{i=1}^{k}
T_i)\Box (\bigcup_{j=1}^{\ell} T_j')}(S)\geq
3\lfloor\frac{k}{2}\rfloor$ since $(\bigcup_{i=1}^{k} T_i)\Box
(\bigcup_{j=1}^{\ell} T_j')$ is a subgraph of $G\Box H$. Therefore,
the problem is converted into finding out
$3\lfloor\frac{k}{2}\rfloor$ internally disjoint pedant $S$-Steiner
trees in $(\bigcup_{i=1}^{k} T_i)\Box (\bigcup_{j=1}^{\ell} T_j')$.
Since
$$
\bigcup_{i=1}^{\lfloor k/2\rfloor} (T_{2i-1}\cup T_{2i})\Box (T_{2i-1}'\cup T_{2i}')
$$
is a subgraph of $(\bigcup_{i=1}^{k} T_i)\Box (\bigcup_{j=1}^{\ell}
T_j')$, we only need to show that
$$
\tau_{G\Box H}(S)\geq
\tau_{\bigcup_{i=1}^{\lfloor k/2\rfloor} (T_{2i-1}\cup T_{2i})\Box (T_{2i-1}'\cup T_{2i}')}(S)\geq 3\lfloor k/2\rfloor.
$$
The structure of $\bigcup_{i=1}^{\lfloor k/2\rfloor} (T_{2i-1}\cup
T_{2i})\Box (T_{2i-1}'\cup T_{2i}')$ in $\bigcup_{i=1}^{\lfloor
k/2\rfloor} (T_{2i-1}\cup T_{2i})\Box H)$ is shown in Fig. 2. In
order to show this structure clearly, we take $2\lfloor k/2\rfloor$
copies of $H(u_j)$, and $2\lfloor k/2\rfloor$ copies of $H(u_k)$.
Note that, these $2\lfloor k/2\rfloor$ copes of $H(u_j)$ (resp.
$H(u_k)$) represent the same graph. For the above example, if we can
prove that $\tau_{(T_1\cup T_2\cup T_3\cup T_4)\Box
\bigcup_{i=1}^{\lfloor k/2\rfloor}(T_{2i-1}'\cup T_{2i}')}(S)\geq
3\lfloor k/2\rfloor$ for $S=\{x,y,z\}$, then $\kappa_{G\Box
H}(S)\geq \tau_{(T_1\cup T_2\cup T_3\cup T_4)\Box
\bigcup_{i=1}^{\lfloor k/2\rfloor}(T_{2i-1}'\cup T_{2i}')}(S)\geq
3\lfloor k/2\rfloor$, as desired. The problem is converted into
finding out $3\lfloor k/2\rfloor$ internally disjoint pedant
$S$-Steiner trees in $(T_1\cup T_2\cup T_3\cup T_4)\Box
\bigcup_{i=1}^{\lfloor k/2\rfloor}(T_{2i-1}'\cup T_{2i}')$ (see Fig.
1 $(h)$).
\begin{figure}[!hbpt]
\begin{center}
\includegraphics[scale=0.8]{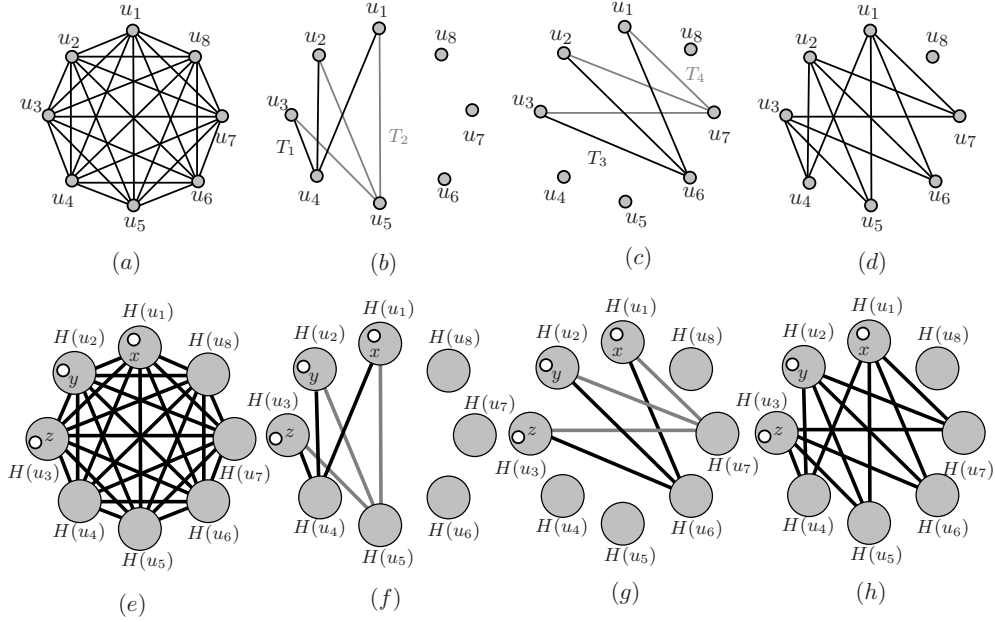}\\
\caption{The structure of $\bigcup_{i=1}^{\lfloor k/2\rfloor} (T_{2i-1}\cup T_{2i})\Box H$.}
\end{center}\label{fig7}
\end{figure}

For each $T_{2i-1}\cup T_{2i}$ and $T_{2i-1}'\cup T_{2i}'$ $(1\leq
i\leq \ell)$, if we can find out $3$ internally disjoint pedant
$S$-Steiner trees in $(T_{2i-1}\cup T_{2i})\Box (T_{2i-1}'\cup
T_{2i}')$, say $T_{i,1},T_{i,2},T_{i,3}$, then the total number of
internally disjoint pedant $S$-Steiner trees in
$\bigcup_{i=1}^{\lfloor k/2\rfloor} (T_{2i-1}\cup T_{2i})\Box
(T_{2i-1}'\cup T_{2i}')$ are $3\lfloor k/2\rfloor$, which implies
that $\tau_{G\Box H}(S)\geq \tau_{\bigcup_{i=1}^{\lfloor k/2\rfloor}
(T_{2i-1}\cup T_{2i})\Box (T_{2i-1}'\cup T_{2i}')}(S)\geq 3\lfloor
k/2\rfloor$ (Note that we must guarantee that any two trees in
$\{T_{i,j}\,|\,1\leq i\leq \lfloor k/2\rfloor, \ 1\leq j\leq 3\}$
are internally disjoint).

Furthermore, from the arbitrariness of $S$, we can get $\tau_3(G\Box
H)\geq 3\lfloor\frac{\tau_3(G)}{2}\rfloor$ and complete the proof of
Theorem \ref{th1-4}. For the above example, we need to find out $3$
internally disjoint pedant $S$-Steiner trees in $(T_{2i-1}\cup
T_{2i})\Box (T_{2i-1}'\cup T_{2i}')$ (see Fig. 1 $(f),(g)$). Then
the total number of internally disjoint  pedant $S$-Steiner in
$\bigcup_{i=1}^{\lfloor k/2\rfloor} (T_{2i-1}\cup T_{2i})\Box
(T_{2i-1}'\cup T_{2i}')$ are $3\lfloor\frac{k}{2}\rfloor$, which
implies $\kappa_{G\Box H}(S)\geq \kappa_{\bigcup_{i=1}^{\lfloor
k/2\rfloor} (T_{2i-1}\cup T_{2i})\Box (T_{2i-1}'\cup
T_{2i}')}(S)\geq 3\lfloor\frac{k}{2}\rfloor$. Thus the result
follows by the arbitrariness of $S$.

From the above analysis, we need to consider the graph $(T_{2i-1}\cup T_{2i})\Box (T_{2i-1}'\cup T_{2i}')$
and prove that for any $S=\{x,y,z\}\subseteq V((T_{2i-1}\cup T_{2i})\Box (T_{2i-1}'\cup T_{2i}'))$ there are
$3$ internally disjoint pedant $S$-Steiner trees in $(T_{2i-1}\cup T_{2i})\Box (T_{2i-1}'\cup T_{2i}')$ for each $i \ (1\leq i\leq \lfloor\frac{k}{2}\rfloor)$.

In the basis of such an idea, we study pedant tree $3$-connectivity
of Cartesian product of the union of two trees $T_1,T_2$ in $G$ and
the union of two trees $T_1',T_2'$ in $H$ first, and show that
$\tau_{3} (T_{2i-1}\cup T_{2i})\Box (T_{2i-1}'\cup T_{2i}')\geq 3$
in Subsection $2.2$. After this preparation, we consider the graph
$G\Box H$ where $G,H$ are two general (connected) graphs and prove
$\kappa_3(G\Box H)\geq \lfloor\frac{\kappa_3(G)}{2}\rfloor$ in
Subsection $2.3$. In Subsection $2.1$, we investigate the pedant
tree $3$-connectivity of Cartesian product of a path $P_n$ and a
connected graph $H$. So the proof of Theorem \ref{th1-4} can be
divided into the above mentioned three subsections. The first and
second subsections are preparations of the last one.

\begin{figure}[!hbpt]
\begin{center}
\includegraphics[scale=0.8]{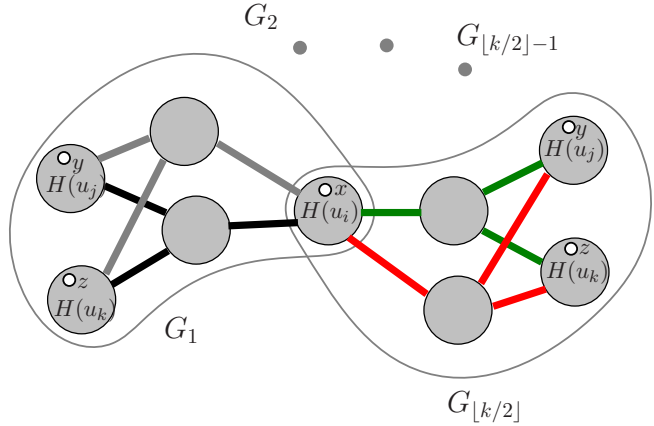}
\caption{Structure of $\bigcup_{i=1}^{\lfloor k/2\rfloor} G_i\Box H$, where $G_i=(T_{2i-1}\cup T_{2i})$.}
\end{center}\label{fig7}
\end{figure}

\subsection{Cartesian product of a path and a connected graph}

A \emph{subdivision} of $G$ is a graph obtained from $G$ by
replacing edges with pairwise internally disjoint paths. If $T$ is
an minimum pedant $S$-Steiner tree, then $T$ is a subdivision of
$K_{1,3}$, and hence $T$ contains a vertex as its root. The
following proposition is a preparation of Subsection 2.3.

\begin{pro}\label{pro2-1}
Let $H$ be a connected graph and $P_n$ be a path with $n$ vertices.
Then $\tau_3(P_n\Box H)\geq \tau_3(H)$. Moreover, the bound is
sharp.
\end{pro}

Set $\tau_3(H)=\ell$, $V(H)=\{v_1,v_2,\ldots,v_m\}$ and
$V(P_n)=\{u_1,u_2,\ldots,u_n\}$. Without loss of generality, let
$u_i$ and $u_j$ be adjacent if and only if $|i-j|=1$, where $1\leq
i\neq j\leq n$. It suffices to show that $\tau_{P_n\Box H}(S)\geq
\ell$ for any $S=\{x,y,z\}\subseteq V(P_n\Box H)$, that is, there
exist $\ell$ internally disjoint pedant $S$-Steiner trees in
$P_n\Box H$. We proceed our proof by the following three lemmas.

\begin{lem}\label{lem2-2}
If $x,y,z$ belongs to the same $V(H(u_i)) \ (1\leq i\leq n)$, then
there exist $\ell+1$ internally disjoint pedant $S$-Steiner trees.
\end{lem}
\begin{pf}
Without loss of generality, we assume $x,y,z\in V(H(u_1))$. Since
$\tau_3(H)=\ell$, it follows that there are $\ell$ internally
disjoint pedant $S$-Steiner trees in $H(u_1)$, say $T_1,T_2,\cdots,
T_{\ell}$. Without loss of generality, let $x=(u_1,v_1)$,
$y=(u_1,v_2)$ and $z=(u_1,v_3)$. Then the trees $T$ induced by the
edges in
$\{x(u_2,v_1),y(u_2,v_2),z(u_2,v_3),(u_2,v_1)(u_2,v_2),(u_2,v_2)(u_2,v_3)\}$
is a pedant $S$-Steiner tree. One can see that the tree $T$ and each
tree $T_i \ (1\leq i\leq \ell)$ are internally disjoint. Therefore,
the trees $T,T_1,T_2,\cdots, T_{\ell}$ are $\ell+1$ internally
disjoint pedant $S$-Steiner trees, as desired.\qed
\end{pf}

\begin{lem}\label{lem2-3}
If only two vertices of $\{x,y,z\}$ belong to some copy $H(u_i) \
(1\leq i\leq n)$, then there exist $\ell$ internally disjoint pedant
$S$-Steiner trees.
\end{lem}
\begin{pf}
We may assume $x,y\in V(H(u_1))$ and $z\in V(H(u_i)) \ (2\leq i\leq
n)$. In the following argument, we can see that this assumption has
no impact on the correctness of our proof. Let $x',y'$ be the
vertices corresponding to $x,y$ in $H(u_i)$, $z'$ be the vertex
corresponding to $z$ in $H(u_1)$.

Suppose $z'\not\in \{x,y\}$. Since $\tau_3(H)=\ell$, it follow that
$\tau_3(H(u_1))=\tau_3(H(u_i))=\ell$, and hence there exist $\ell$
internally disjoint pedant $S$-Steiner trees
$T_1,T_2,\cdots,T_{\ell}$ in $H(u_1)$ and there exist $\ell$
internally disjoint pedant $S$-Steiner trees
$T_1',T_2',\cdots,T_{\ell}'$ in $H(u_i)$. Let $w_i,w_i'$ be the root
of $T_i,T_i'$, respectively. Let $P_i,Q_i,R_i$ denote the unique
path connecting $w_i$ and $x,y,z'$, respectively. Let
$P_i',Q_i',R_i'$ denote the unique path connecting $w_i'$ and
$x',y',z$, respectively. Without loss of generality, let
$w_i=(u_1,v_1)$ and $w_i'=(u_i,v_1)$. Then the trees $T_i$ induced
by the edges in $E(P_i)\cup E(Q_i)\cup E(R_i')\cup
\{(u_j,v_1)(u_{j+1},v_1)\,|\,1\leq j\leq i-1\} \ (1\leq i\leq \ell)$
are $\ell$ internally disjoint pedant $S$-Steiner trees.

\begin{figure}[!hbpt]
\begin{center}
\includegraphics[scale=1]{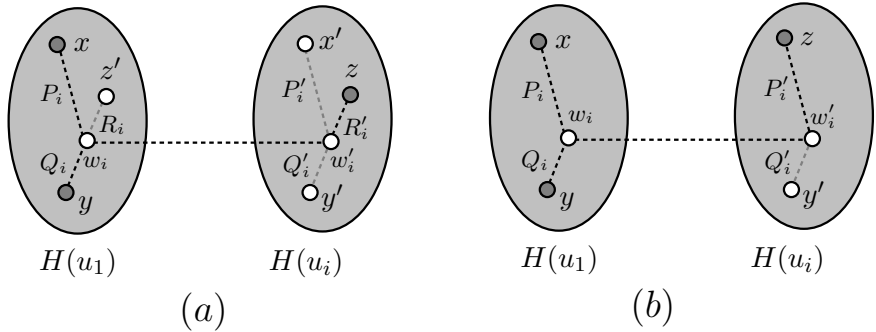}
\end{center}
\begin{center}
\caption{Graphs for Lemma \ref{lem2-3}.}
\end{center}\label{fig7}
\end{figure}

Suppose $z'\in \{x,y\}$. Without loss of generality, let $z'=x$.
Since $\tau_3(H)=\ell$, it follow from Lemma \ref{lem1-2} that
$\kappa(H)\geq \ell+1$, and hence $\kappa(H(u_1))\geq \ell+1$ and
$\kappa(H(u_i))\geq \ell+1$. Then there exist $\ell+1$ internally
disjoint paths connecting $x$ and $y$ in $H(u_1)$, say
$R_1,R_2,\cdots,R_{\ell+1}$, and there exist $\ell$ internally
disjoint paths connecting $z$ and $y'$ in $H(u_i)$, say
$R_1',R_2',\cdots,R_{\ell+1}'$. Note that there is at most one path
in $\{R_1,R_2,\cdots,R_{\ell+1}\}$, say $R_{\ell+1}$, such that its
length is $1$, and there is at most one path in
$\{R_1',R_2',\cdots,R_{\ell+1}'\}$, say $R_{\ell+1}'$, such that its
length is $1$. Then there is an internal vertex $w_i$ in $R_i$, and
there is an internal vertex $w_i'$ in $R_i'$. Let $P_i,Q_i$ denote
the unique path connecting $w_i$ and $x,y$, respectively. Let
$P_i',Q_i'$ denote the unique path connecting $w_i'$ and $y',z$,
respectively. Without loss of generality, let $w_i=(u_1,v_1)$ and
$w_i'=(u_i,v_1)$. Then the trees $T_i$ induced by the edges in
$E(P_i)\cup E(Q_i)\cup E(P_i')\cup
\{(u_j,v_1)(u_{j+1},v_1)\,|\,1\leq j\leq i-1\} \ (1\leq i\leq \ell)$
are $\ell$ internally disjoint pedant $S$-Steiner trees, as
desired.\qed
\end{pf}

\begin{lem}\label{lem2-4}
If $x,y,z$ are contained in distinct $H(u_i)$s, then there exist
$\ell$ internally disjoint pedant $S$-Steiner trees.
\end{lem}
\begin{pf}
We may assume that $x\in V(H(u_1))$, $y\in V(H(u_i))$, $z\in
V(H(u_j))$. In the following argument, we can see that this
assumption has no influence on the correctness of our proof. Let
$y',z'$ be the vertices corresponding to $y,z$ in $H(u_1)$, $x',z''$
be the vertices corresponding to $x,z$ in $H(u_i)$ and $x'',y''$ be
the vertices corresponding to $x,y$ in $H(u_j)$.

Suppose that $x,y',z'$ are distinct vertices in $H(u_1)$. Since
$\tau_3(H)=\ell$, it follow that
$\tau_3(H(u_1))=\tau_3(H(u_i))=\tau_3(H(u_j))=\ell$, and hence there
exist $\ell$ internally disjoint pedant $S$-Steiner trees
$T_1,T_2,\cdots,T_{\ell}$ in $H(u_1)$, and there exist $\ell$
internally disjoint pedant $S$-Steiner trees
$T_1',T_2',\cdots,T_{\ell}'$ in $H(u_i)$, and there exist $\ell$
internally disjoint pedant $S$-Steiner trees
$T_1'',T_2'',\cdots,T_{\ell}''$ in $H(u_j)$. Let $w_i,w_i',w_i''$ be
the root of $T_i,T_i',T_i''$, respectively. Let $P_i,Q_i,R_i$ denote
the unique path connecting $w_i$ and $x,y,z'$, respectively. Let
$P_i',Q_i',R_i'$ denote the unique path connecting $w_i'$ and
$x,y',z'$, respectively. Let $P_i'',Q_i'',R_i''$ denote the unique
path connecting $w_i''$ and $x',y,z''$, respectively. Without loss
of generality, let $w_i=(u_1,v_1)$, $w_i'=(u_i,v_1)$ and
$w_i''=(u_j,v_1)$. Then the trees $T_i$ induced by the edges in
$E(P_i)\cup E(Q_i')\cup E(R_i'')\cup
\{(u_r,v_1)(u_{r+1},v_1)\,|\,1\leq r\leq i-1\}\cup
\{(u_r,v_1)(u_{r+1},v_1)\,|\,i\leq r\leq j-1\} \ (1\leq i\leq \ell)$
are $\ell$ internally disjoint pedant $S$-Steiner trees.
\begin{figure}[!hbpt]
\begin{center}
\includegraphics[scale=1]{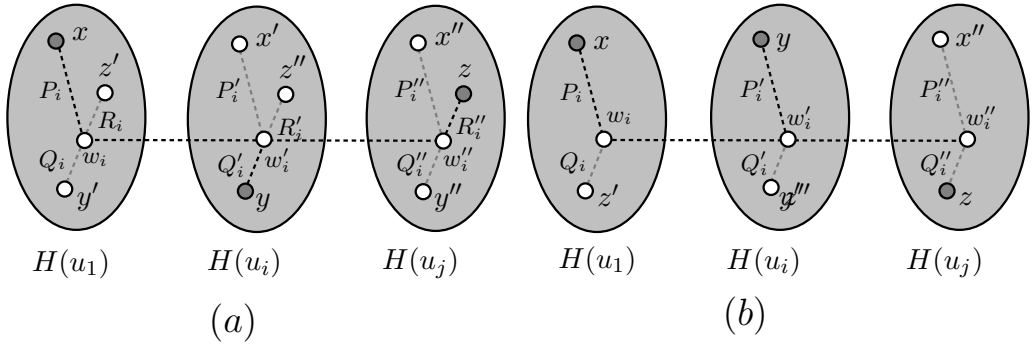}
\end{center}
\begin{center}
\caption{Graphs for Case $1$ of Lemma \ref{lem2-4}}.
\end{center}\label{fig7}
\end{figure}

Suppose that two of $x, y',z'$ are the same vertex in $H(u_1)$.
Without loss of generality, let $x=y'$. Since $\tau_3(H)=\ell$, it
follow from Lemma \ref{lem1-2} that $\kappa(H)\geq \ell+1$, and
hence $\kappa(H(u_{1}))\geq \ell+1$, $\kappa(H(u_i))\geq \ell+1$ and
$\kappa(H(u_j))\geq \ell+1$. Then there exist $\ell+1$ internally
disjoint paths connecting $x$ and $z'$ in $H(u_1)$, say
$R_1,R_2,\cdots,R_{\ell+1}$, and there exist $\ell$ internally
disjoint paths connecting $y$ and $z''$ in $H(u_i)$, say
$R_1',R_2',\cdots,R_{\ell+1}'$, and there exist $\ell$ internally
disjoint paths connecting $y$ and $z''$ in $H(u_j)$, say
$R_1'',R_2'',\cdots,R_{\ell+1}''$. Note that there is at most one
path in $\{R_1,R_2,\cdots,R_{\ell+1}\}$, say $R_{\ell+1}$, such that
its length is $1$, and there is at most one path in
$\{R_1',R_2',\cdots,R_{\ell+1}'\}$, say $R_{\ell+1}'$, such that its
length is $1$, and there is at most one path in
$\{R_1'',R_2'',\cdots,R_{\ell+1}''\}$, say $R_{\ell+1}''$, such that
its length is $1$. Then there is an internal vertex $w_i$ in $R_i$,
and there is an internal vertex $w_i'$ in $R_i'$, and there is an
internal vertex $w_i''$ in $R_i''$. Let $P_i,Q_i$ denote the unique
path connecting $w_i$ and $x,z'$, respectively. Let $P_i',Q_i'$
denote the unique path connecting $w_i'$ and $y,z''$, respectively.
Let $P_i'',Q_i''$ denote the unique path connecting $w_i'$ and
$x'',z$, respectively. Then the trees $T_i$ induced by the edges in
$E(P_i)\cup E(P_i')\cup E(Q_i'')\cup
\{(u_j,v_1)(u_{j+1},v_1)\,|\,1\leq j\leq i-1\}\cup
\{(u_j,v_1)(u_{j+1},v_1)\,|\,i\leq j\leq j-1\} \ (1\leq i\leq \ell)$
are $\ell$ internally disjoint pedant $S$-Steiner trees, as desired.

Suppose that $x,y',z'$ are the same vertex in $H(u_1)$.
Since $\tau_3(H)=\ell$, it follow from Lemma \ref{lem1-1} that $\delta(H)\geq \ell+2$, and hence $\delta(H(u_1))\geq \ell+1$, $\delta(H(u_i))\geq \ell+1$ and $\delta(H(u_j))\geq \ell+1$. Then there are $\ell+1$ neighbors of $x$ in $H(u_1)$, say $(u_1,v_1),(u_1,v_2),\cdots,(u_1,v_{\ell+1})$. By the same reason, there are $\ell+1$ neighbors of $y$ in $H(u_i)$, say $(u_i,v_1),(u_i,v_2),\cdots,(u_i,v_{\ell+1})$, and there are $\ell+1$ neighbors of $z$ in $H(u_j)$, say $(u_j,v_1),(u_j,v_2),\cdots,(u_j,v_{\ell+1})$. Then the tree $T_r$ induced by the edges in $\{x(u_1,v_r),y(u_i,v_r),z(u_j,v_r)\}\cup \{(u_{s},v_r)(u_{s+1},v_r)\,|\,1\leq s\leq i-1\}\cup \{(u_{s},v_r)(u_{s+1},v_r)\,|\,i\leq s\leq j-1\}$ is a pedant $S$-Steiner tree, where $1\leq r\leq \ell+1$. Therefore, the trees $T_1,T_2,\cdots,T_{\ell+1}$ are
$\ell+1$ internally disjoint pedant $S$-Steiner trees, as desired. \qed
\end{pf}

From Lemmas \ref{lem2-2}, \ref{lem2-3} and \ref{lem2-4}, we conclude
that, for any $S\subseteq V(P_n\Box H)$, there exist $\ell$
internally disjoint pendant $S$-Steiner trees, which implies that
$\tau_{P_n\Box H}(S)\geq \ell$. From the arbitrariness of $S$, we
have $\tau_3(P_n\Box H)\geq \ell$.
The proof of Proposition \ref{pro2-1} is complete.

\subsection{Cartesian product of two trees in $G$ and two trees in $H$}

In this subsection, we consider the pedant tree $3$-connectivity of
Cartesian product of two trees in $G$ and two trees in $H$, which is
a preparation of the next subsection.

\begin{pro}\label{pro2-5}
Let $G,H$ be two graphs. For $S=\{x,y,z\}\subseteq V(G\Box H)$, we
assume that $x\in V(H(u_1))$, $y\in V(H(u_2))$ and $z\in V(H(u_3))$.
Let $T_1,T_2$ be two minimum pedant Steiner trees connecting
$\{u_1,u_2,u_3\}$ in $G$. Let $y',z'$ be the vertices corresponding
to $y,z$ in $H(u_1)$. Let $T_1',T_2'$ be two pedant Steiner trees
connecting $\{x,y',z'\}$ in $H$. Then
$$
\tau_{(T_1\cup T_2)\Box (T_1'\cup T_2')}(S)\geq 3.
$$
\end{pro}
\begin{pf}
Since $T_1,T_2$ are two minimum pedant Steiner trees connecting
$\{u_1,u_2,u_3\}$, it follows that $T_1,T_2$ are subdivisions of
$K_{1,3}$ and hence have roots, say $u_s,u_t$, respectively. Note
that $x\in V(H(u_{1}))$, $y\in V(H(u_2))$ and $z\in V(H(u_3))$. Let
$y',z'$ be the vertices corresponding to $y,z$ in $H(u_1)$, $x',z''$
be the vertices corresponding to $x,z$ in $H(u_2)$ and $x'',y''$ be
the vertices corresponding to $x,y$ in $H(u_3)$.

$\bullet$ Let $R_1,R_2,R_3$ be the three paths connecting $u_r$ and $u_1,u_2,u_3$, respectively.

~~~~~~$\bullet$ Set $R_1=u_1p_1p_2\cdots p_au_r$, where $p_i\in V(G), 1\leq i\leq a$.

~~~~~~$\bullet$ Set $R_2=u_2p_1'p_2'\cdots p_b'u_r$, where $p_i'\in V(G), 1\leq i\leq b$.

~~~~~~$\bullet$ Set $R_3=u_3p_1''p_2''\cdots p_c''u_r$, where $p_i''\in V(G), 1\leq i\leq c$.

$\bullet$ Let $R_1',R_2',R_3'$ be the three paths connecting $u_s$ and $u_1,u_2,u_3$, respectively.

~~~~~~$\bullet$ Set $R_1'=u_1q_1q_2\cdots q_du_r$, where $q_i\in V(G), 1\leq i\leq d$.

~~~~~~$\bullet$ Set $R_2'=u_2q_1'q_2'\cdots q_e'u_r$, where $q_i'\in V(G), 1\leq i\leq e$.

~~~~~~$\bullet$ Set $R_3'=u_3q_1''q_2''\cdots q_f''u_r$, where $q_i''\in V(G), 1\leq i\leq f$.

We distinguish the following three cases to show this proposition.

{\bf Case 1.} The vertices $x,y',z'$ are distinct vertices in $H(u_1)$.

In order to show the structure of pedant $S$-Steiner trees clearly,
we assume all of the following.

$\bullet$ Let $s,t$ be the roots of $T_1',T_2'$, respectively.

$\bullet$ Let $s',s'',s_1,s_2$ be the vertices corresponding to $s$
in $H(u_2),H(u_3),H(u_r),H(u_s)$, respectively.

$\bullet$ Let $t',t'',t_1,t_2$ be the vertices corresponding to $t$
in $H(u_2),H(u_3),H(u_r),H(u_s)$, respectively.

$\bullet$ Let $P_{1,1},P_{1,2},P_{1,3}$ be the three paths
connecting $s$ and $x,y',z'$ in $T_1'$, respectively.

$\bullet$ Let $Q_{1,1},Q_{1,2},Q_{1,3}$ be the three paths
connecting $s$ and $x,y',z'$ in $T_1'$, respectively.

$\bullet$ Let $P_{2,j},P_{3,j},P_{r,j},P_{s,j} \ (1\leq j\leq 3)$ be
the paths corresponding to $P_{1,j}$ in $H(u_2),H(u_3),$
$H(u_r),H(u_s)$, respectively.

$\bullet$ Let $Q_{2,j},Q_{3,j},Q_{r,j},Q_{s,j} \ (1\leq j\leq 3)$ be
the paths corresponding to $Q_{1,j}$ in $H(u_2),H(u_3),$
$H(u_r),H(u_s)$, respectively.

$\bullet$ Without loss of generality, let $x=(u_1,v_1)$,
$y'=(u_1,v_2)$, $z'=(u_1,v_3)$, $s=(u_1,v_4)$ and $t=(u_1,v_5)$.

\begin{figure}[!hbpt]
\begin{center}
\includegraphics[scale=0.9]{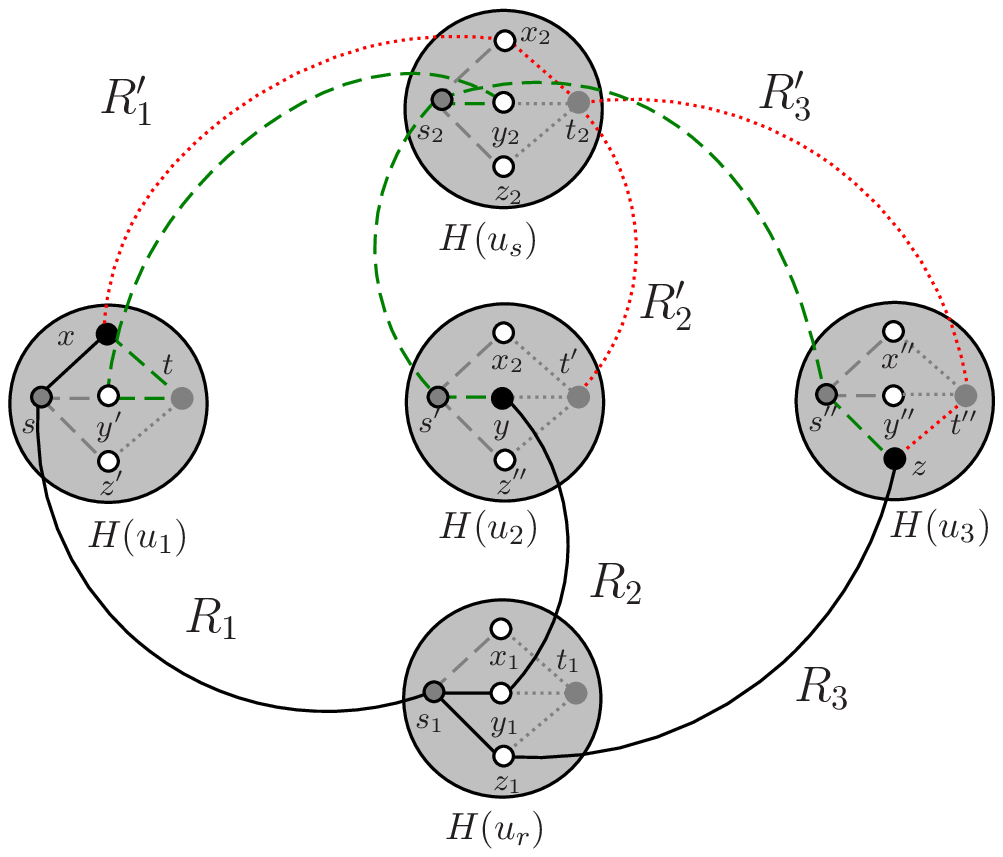}
\end{center}
\begin{center}
\caption{Graphs for Case 1 of Proposition \ref{pro2-5}.}
\end{center}\label{fig7}
\end{figure}

Then the tree $T$ induced by the edges in
\begin{eqnarray*}
&&E(P_{1,1})\cup E(P_{r,2})\cup E(P_{r,3})\\
&&\cup \{s(p_1,v_4)\}\cup \{(p_i,v_4)(p_{i+1},v_4)\,|\,1\leq i\leq a-1\}\cup \{(p_a,v_4)s_1\}\\
&&\cup \{y(p_1',v_2)\}\cup \{(p_i',v_2)(p_{i+1}',v_2)\,|\,1\leq i\leq b-1\}\cup \{(p_b',v_2)y_1\}\\
&&\cup \{z(p_1'',v_3)\}\cup \{(p_i'',v_3)(p_{i+1}'',v_3)\,|\,1\leq i\leq c-1\}\cup \{(p_c'',v_3)z_1\}
\end{eqnarray*}
is a pedant $S$-Steiner tree, and the tree $T'$ induced by the edges in
\begin{eqnarray*}
&&E(Q_{s,1})\cup E(Q_{2,2})\cup E(Q_{3,3})\\
&&\cup \{x(q_1,v_1)\}\cup \{(q_i,v_1)(q_{i+1},v_1)\,|\,1\leq i\leq d-1\}\cup \{(q_d,v_1)x_2\}\\
&&\cup \{t'(q_1',v_5)\}\cup \{(q_i',v_5)(q_{i+1}',v_5)\,|\,1\leq i\leq e-1\}\cup \{(q_e',v_5)t_2\}\\
&&\cup \{t''(q_1'',v_5)\}\cup \{(q_i'',v_5)(q_{i+1}'',v_5)\,|\,1\leq
i\leq f-1\}\cup \{(q_c'',v_5)t_2\}
\end{eqnarray*}
is a pedant $S$-Steiner tree, and the tree $T''$ induced by the edges in
\begin{eqnarray*}
&&E(Q_{s,2})\cup E(Q_{1,2})\cup E(Q_{1,1})\cup E(P_{2,2})\cup E(P_{3,3})\\
&&\cup \{y'(q_1,v_2)\}\cup \{(q_i,v_2)(q_{i+1},v_2)\,|\,1\leq i\leq d-1\}\cup \{(q_d,v_2)y_2\}\\
&&\cup \{s'(q_1',v_4)\}\cup \{(q_i',v_4)(q_{i+1}',v_4)\,|\,1\leq i\leq e-1\}\cup \{(q_e',v_4)s_2\}\\
&&\cup \{s''(q_1'',v_4)\}\cup \{(q_i'',v_4)(q_{i+1}'',v_4)\,|\,1\leq
i\leq f-1\}\cup \{(q_f'',v_4)s_2\}
\end{eqnarray*}
is a pedant $S$-Steiner tree. Since $T,T',T''$ are internally disjoint, it follows that
$$
\tau_{(T_1\cup T_2)\Box (T_1'\cup T_2')}(S)\geq 3,
$$
as desired.

{\bf Case 2.} Two of $x, y',z'$ are the same vertex in $H(u_1)$.

Without loss of generality, let $x=z'$. Note that there are two
paths $P_1,Q_1$ connecting $x$ and $y'$ in $T_1,T_2$, respectively.
Observe that at most of $P_1,Q_1$ is length $1$. We now assume that
the length of $P_1$ is at least $2$. Then there exists an internal
vertex in $P_1$, say $s$, and hence $s$ divide $P_1$ into two paths,
say $P_{1,1},P_{1,2}$. In order to show the structure of pedant
$S$-Steiner tree clearly, we assume the following.

$\bullet$ Let $x_1,x_2$ be the vertices corresponding to $x$ in $H(u_r),H(u_s)$, respectively.

$\bullet$ Let $y_1,y_2$ be the vertices corresponding to $y'$ in
$H(u_r),H(u_s)$, respectively.

$\bullet$ Let $s',s'',s_1,s_2$ be the vertices corresponding to $s$
in $H(u_2),H(u_3),H(u_r),H(u_s)$, respectively.

$\bullet$ Without loss of generality, let $x=(u_1,v_1)$, $y'=(u_1,v_2)$ and $s=(u_1,v_3)$.

$\bullet$ Let $P_{2,j},P_{3,j},P_{r,j},P_{s,j} \ (j=1,2)$ be the
paths corresponding to $P_{1,j}$ in $H(u_2),H(u_3),$
$H(u_r),H(u_s)$, respectively.

$\bullet$ Let $Q_{2},Q_{3},Q_{r},Q_{s}$ be the paths corresponding
to $Q_{1}$ in $H(u_2),H(u_3),H(u_r),H(u_s)$, respectively.
\begin{figure}[!hbpt]
\begin{center}
\includegraphics[scale=0.9]{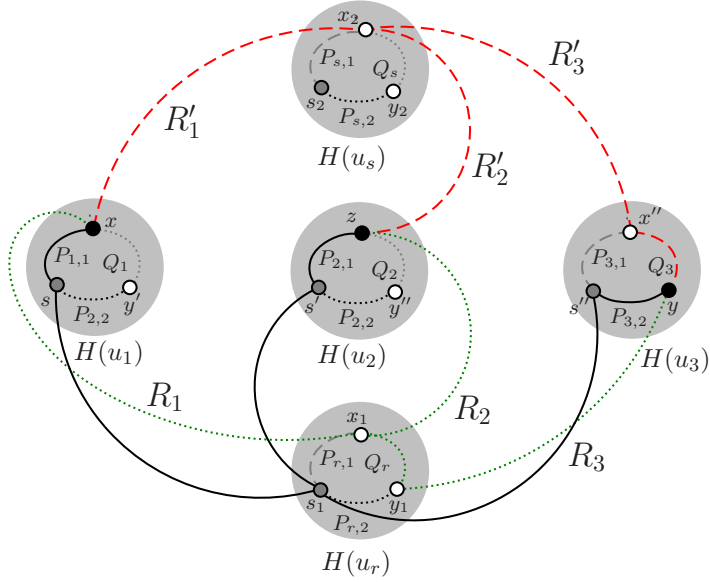}
\end{center}
\begin{center}
\caption{Graphs for Case 2 of Proposition \ref{pro2-5}.}
\end{center}\label{fig7}
\end{figure}

Then the tree $T$ induced by the edges in
\begin{eqnarray*}
&&\cup \{x(p_1,v_1)\}\cup \{(p_i,v_1)(p_{i+1},v_1)\,|\,1\leq i\leq a-1\}\cup \{(p_a,v_1)x_1\}\\
&&\cup \{z(p_1',v_1)\}\cup \{(p_i',v_1)(p_{i+1}',v_1)\,|\,1\leq i\leq b-1\}\cup \{(p_b',v_1)x_1\}\\
&&\cup E(Q_r)\cup \{y(p_1'',v_2)\}\cup \{(p_i'',v_2)(p_{i+1}'',v_2)\,|\,1\leq i\leq c-1\}\cup \{(p_c'',v_2)y_1\}
\end{eqnarray*}
is a pedant $S$-Steiner tree, and the tree $T'$ induced by the edges in
\begin{eqnarray*}
&&\cup \{x(q_1,v_1)\}\cup \{(q_i,v_1)(q_{i+1},v_1)\,|\,1\leq i\leq d-1\}\cup \{(q_d,v_1)x_2\}\\
&&\cup \{z(q_1',v_1)\}\cup \{(q_i',v_1)(q_{i+1}',v_1)\,|\,1\leq i\leq e-1\}\cup \{(q_e',v_1)x_2\}\\
&&\cup \{x''(q_1'',v_1)\}\cup \{(q_i'',v_1)(q_{i+1}'',v_1)\,|\,1\leq i\leq f-1\}\cup \{(q_f'',v_1)x_2\}\cup E(Q_3)
\end{eqnarray*}
is a pedant $S$-Steiner tree, and the tree $T''$ induced by the edges in
\begin{eqnarray*}
&&E(P_{1,1})\cup E(P_{2,1})\cup E(P_{3,2})\\
&&\cup \{s(p_1,v_3)\}\cup \{(p_i,v_3)(p_{i+1},v_3)\,|\,1\leq i\leq a-1\}\cup \{(p_a,v_3)s_1\}\\
&&\cup \{s'(p_1',v_3)\}\cup \{(p_i',v_3)(p_{i+1}',v_3)\,|\,1\leq i\leq b-1\}\cup \{(p_b',v_3)s_1\}\\
&&\cup \{s''(p_1'',v_3)\}\cup \{(p_i'',v_3)(p_{i+1}'',v_3)\,|\,1\leq i\leq c-1\}\cup \{(p_c'',v_3)s_1\}
\end{eqnarray*}
is a pedant $S$-Steiner tree. Since $T,T',T''$ are internally disjoint, it follows that
$$
\tau_{(T_1\cup T_2)\Box (T_1'\cup T_2')}(S)\geq 3,
$$
as desired.

{\bf Case 3.} $x,y',z'$ are the same vertex in $H(u_1)$.

In this section, we let $s,t$ be the neighbors of $x$ in $T_1',T_2'$, respectively. Let $s',s'',s_1,s_2$ be the vertices corresponding to $s$ in $H(u_2),H(u_3),H(u_r),H(u_s)$, respectively. Let $t',t'',t_1,t_2$ be the vertices corresponding to $t$ in $H(u_2),H(u_3),H(u_r),H(u_s)$, respectively.
Without loss of generality, let $x=(u_1,v_1)$, $s=(u_1,v_2)$ and $t=(u_1,v_3)$.
\begin{figure}[!hbpt]
\begin{center}
\includegraphics[scale=0.9]{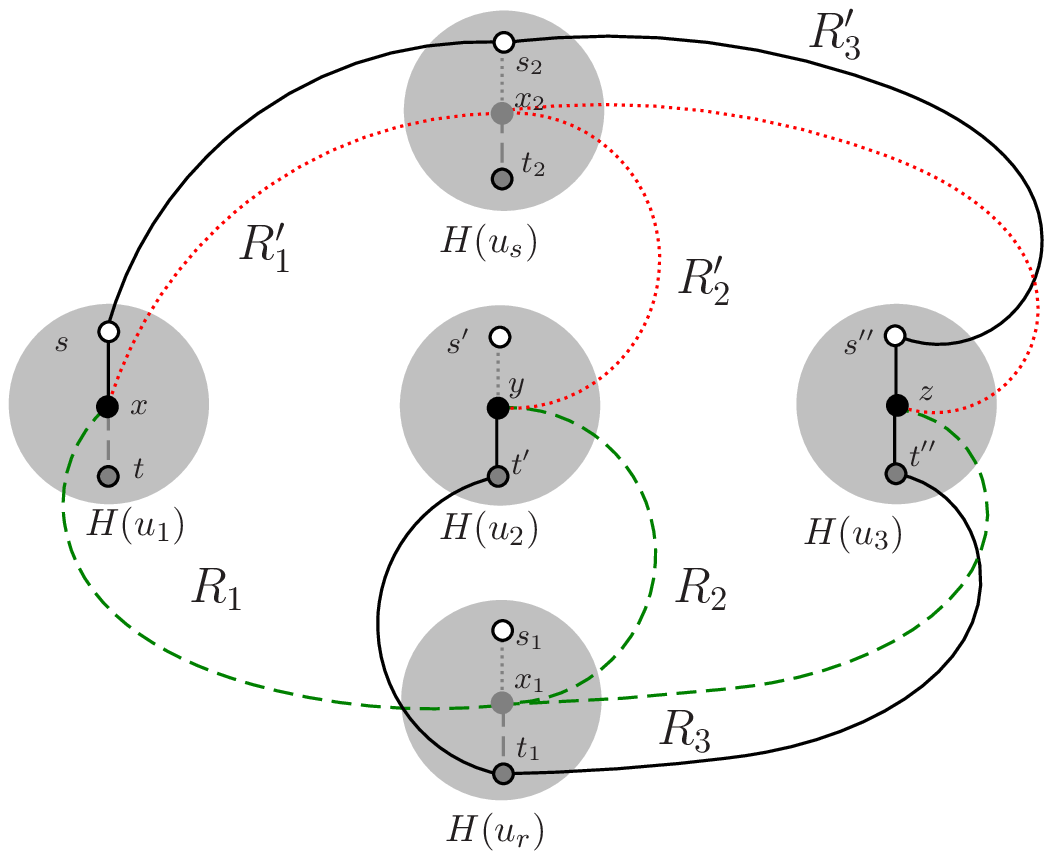}
\end{center}
\begin{center}
\caption{Graphs for Case 3 of Proposition \ref{pro2-5}.}
\end{center}\label{fig7}
\end{figure}

Then the tree $T$ induced by the edges in
\begin{eqnarray*}
&&\cup \{x(p_1,v_1)\}\cup \{(p_i,v_1)(p_{i+1},v_1)\,|\,1\leq i\leq a-1\}\cup \{(p_a,v_1)x_1\}\\
&&\cup \{y(p_1',v_1)\}\cup \{(p_i',v_1)(p_{i+1}',v_1)\,|\,1\leq i\leq b-1\}\cup \{(p_b',v_1)x_1\}\\
&&\cup \{z(p_1'',v_1)\}\cup \{(p_i'',v_1)(p_{i+1}'',v_1)\,|\,1\leq i\leq c-1\}\cup \{(p_c'',v_1)x_1\}
\end{eqnarray*}
is a pedant $S$-Steiner tree, and the tree $T'$ induced by the edges in
\begin{eqnarray*}
&&\cup \{x(q_1,v_1)\}\cup \{(q_i,v_1)(q_{i+1},v_1)\,|\,1\leq i\leq d-1\}\cup \{(q_d,v_1)x_2\}\\
&&\cup \{y(q_1',v_1)\}\cup \{(q_i',v_1)(q_{i+1}',v_1)\,|\,1\leq i\leq e-1\}\cup \{(q_e',v_1)x_2\}\\
&&\cup \{z(q_1'',v_1)\}\cup \{(q_i'',v_1)(q_{i+1}'',v_1)\,|\,1\leq i\leq f-1\}\cup \{(q_f'',v_1)x_2\}
\end{eqnarray*}
is a pedant $S$-Steiner tree, and the tree $T''$ induced by the edges in
\begin{eqnarray*}
&&\{xs,s(q_1,v_2)\}\cup \{(q_i,v_2)(q_{i+1},v_2)\,|\,1\leq i\leq d-1\}\cup \{(q_d,v_2)s_2\}\\
&&\cup \{s_2(q_1'',v_2)\}\cup \{(q_i'',v_2)(q_{i+1}'',v_2)\,|\,1\leq i\leq f-1\}\cup \{(q_f'',v_2)s''\}\\
&&\cup \{s''z,zt'',t''(p_1'',v_3)\}\cup \{(p_i'',v_3)(p_{i+1}'',v_3)\,|\,1\leq i\leq c-1\}\cup \{(p_c'',v_3)t_1\}\\
&&\cup \{yt',t'(p_1',v_3)\}\cup \{(p_i',v_3)(p_{i+1}',v_3)\,|\,1\leq i\leq b-1\}\cup \{(p_b'',v_3)t_1\}
\end{eqnarray*}
is a pedant $S$-Steiner tree. Since $T,T',T''$ are internally disjoint, we have
$$
\tau_{(T_1\cup T_2)\Box (T_1'\cup T_2')}(S)\geq 3,
$$
as desired.

From the above argument,
there exist $3$ internally disjoint pedant $S$-Steiner trees,
which implies $\tau_{T\Box H}(S)\geq 3$.
The proof is now complete.\qed
\end{pf}

\subsection{Cartesian product of two general graphs}

After the above preparations, we are ready to prove Theorem
\ref{th1-4} in this subsection.\\

\noindent\textbf{Proof of Theorem \ref{th1-4}:} Set $\tau_3(G)=k$
and $\tau_3(H)=\ell$. Without loss of generality, let $k\leq \ell$.
Recall that $V(G)=\{u_1,u_2,\ldots,u_n\}$,
$V(H)=\{v_1,v_2,\ldots,v_m\}$. From the definition of $\tau_3(G\Box
H)$ and the symmetry of Cartesian product graphs, we need to prove
that $\tau_{G\Box H}(S)\geq 3\lfloor k/2\rfloor$ for any
$S=\{x,y,z\}\subseteq V(G\Box H)$. Furthermore, it suffices to show
that there exist $3\lfloor k/2\rfloor$ internally disjoint pedant
$S$-Steiner trees in $G\Box H$. Clearly, $V(G\Box
H)=\bigcup_{i=1}^nV(H(u_i))$. Without loss of generality, let $x\in
V(H(u_i))$, $y\in V(H(u_j))$ and $z\in V(H(u_k))$.

{\bf Case $1$}. The vertices $x,y,z$ belongs to the same $V(H(u_i)) \ (1\leq i\leq
n)$.

Without loss of generality, let $x,y,z\in V(H(u_1))$. From
Lemma \ref{lem1-1}, $\delta(G)\geq \tau_3(G)+2=k+2$ and hence the
vertex $u_1$ has at least $k+2$ neighbors in $G$. Select $k+2$ neighbors from
them, say $u_2,u_3,\cdots,u_{k+3}$. Without loss of generality, let
$x=(u_1,v_1)$, $y=(u_1,v_2)$ and $z=(u_1,v_3)$. Clearly, the trees
$T_{i}$ induced by the edges in $\{x(u_i,v_1),
y(u_i,v_2),z(u_i,v_3),(u_i,v_1)$ $(u_i,v_2),(u_i,v_2)(u_i,v_3)\} \
(2\leq i\leq k+3)$ are $k+2$ internally disjoint pedant $S$-Steiner
trees in $G\Box H$, which occupy no edge of $H(u_1)$. Since
$\tau_3(H)=\ell$, it follows that there are $\ell$ internally disjoint
pedant $S$-Steiner trees in $H(u_1)$. Observe that these $\ell$
pedant $S$-Steiner trees and the trees $T_{i} \ (2\leq i\leq k+3)$
are internally disjoint. So the total number of internally disjoint
pedant $S$-Steiner trees is $k+\ell+2>3\lfloor k/2\rfloor$, as desired.

{\bf Case $2$}. Only two vertices of $\{x,y,z\}$ belong to some copy
$H(u_j) \ (1\leq j\leq n)$.

Without loss of generality, let $x,y\in H(u_1)$ and $z\in H(u_2)$.
From Lemma \ref{lem1-2}, $\kappa(G)\geq \tau_3(G)+1=k+1$ and hence
there exist $k+1$ internally disjoint paths connecting $u_1$ and
$u_2$ in $G$, say $P_1,P_2,\cdots,P_{k+1}$. Clearly, there exists at
most one of $P_1,P_2,\cdots,P_{k+1}$, say $P_{k+1}$, such that
$P_{k+1}=u_1u_2$. We may assume that the length of $P_i \ (1\leq
i\leq \ell)$ is at least $2$. From Proposition \ref{pro2-1}, there
exist $\ell$ internally disjoint pedant $S$-Steiner trees in
$P_{k+1}\Box H$, say $T_1,T_2,\cdots,T_{\ell}$. For each $P_i \
(1\leq i\leq \ell)$, since $P_i$ is a path of length at least $2$,
it follows that there exists an internal vertex in $P_i$, say $u_i$.
Let $Q_i,R_i$ be the two paths connecting $u_i$ and $u_1,u_2$ in
$P_i$, respectively. Set $Q_i=u_1,u_1',u_2',\cdots,u_s',u_i$ and
$R_i=u_2,u_1'',u_2'',\cdots,u_t'',u_i$. In the following argument,
we can see that this assumption has no impact on the correctness of
our proof. Let $x',y'$ be the vertices corresponding to $x,y$ in
$H(u_2)$, $z'$ be the vertex corresponding to $z$ in $H(u_1)$, and
$x'',y'',z''$ be the vertices corresponding to $x,y,z$ in $H(u_i)$.

Suppose $z'\not\in \{x,y\}$. Without loss of generality, let
$x=(u_1,v_1)$, $y=(u_1,v_2)$ and $z=(u_2,v_3)$. Since $\tau_3(G)\geq
\ell\geq 1$, it follows that there is a pedant Steiner tree
connecting $\{x'',y'',z''\}$ in $H(u_i)$, say $T^{i}$. Furthermore,
the tree $T_i' \ (1\leq i\leq k)$ induced by the edges in
\begin{eqnarray*}
&&E(T^{i})\cup \{x(u_1',v_1)\}\cup \{(u_j',v_1)(u_{j+1}',v_1)\,|\,1\leq i\leq s\}\cup \{x''(u_s',v_1)\}\cup \{y(u_1',v_2)\}\\
&&\cup \{(u_j',v_2)(u_{j+1}',v_2)\,|\,1\leq i\leq s\}\cup \{x''(u_s',v_2)\}\cup \{z(u_1'',v_3)\}\\
&&\cup \{(u_j'',v_2)(u_{j+1}'',v_2)\,|\,1\leq i\leq t\}\cup \{z''(u_s'',v_3)\}\\
\end{eqnarray*}
is a pedant $S$-Steiner tree. Obviously, the trees
$T_1,T_2,\cdots,T_{\ell},T_1',T_2',\cdots,T_{k}'$ are $k+\ell\geq
3\lfloor k/2 \rfloor$ internally disjoint pedant $S$-Steiner trees.

Suppose $z'\in \{x,y\}$. Without loss of generality, let $z'=x$.
$x=(u_1,v_1)$, $y=(u_1,v_2)$. Then $z=(u_2,v_1)$. Since
$\tau_3(G)\geq 1$, it follows that there is a path connecting $x''$
and $y''$, say $P'$. Furthermore, the tree $T_i' \ (1\leq i\leq k)$
induced by the edges in
\begin{eqnarray*}
&&E(P')\cup \{x(u_1',v_1)\}\cup \{(u_j',v_1)(u_{j+1}',v_1)\,|\,1\leq i\leq s\}\cup \{x''(u_s',v_1)\}\cup \{y(u_1',v_2)\}\\
&&\cup \{(u_j',v_2)(u_{j+1}',v_2)\,|\,1\leq i\leq s\}\cup \{x''(u_s',v_2)\}\cup \{z(u_1'',v_3)\}\\
&&\cup \{(u_j'',v_2)(u_{j+1}'',v_2)\,|\,1\leq i\leq t\}\cup \{z''(u_s'',v_3)\}\\
\end{eqnarray*}
is a pedant $S$-Steiner tree. Obviously, the trees $T_1,T_2,\cdots,T_{\ell},T_1',T_2',\cdots,T_{k}'$ are $k+\ell\geq 3\lfloor k/2 \rfloor$ internally disjoint pedant $S$-Steiner trees.

{\bf Case $3$}. The vertices $x,y,z$ are contained in distinct $H(u_i)$s.

Without loss of generality, let $x\in V(H(u_1))$, $y\in V(H(u_2))$
and $z\in V(H(u_3))$. Since $\tau_3(G)=k$, it follows that there
exist $k$ internally disjoint pedant Steiner trees connecting
$\{u_1,u_2,u_3\}$ in $G$, say $T_1,T_2,\cdots,T_{\ell}$. Let $y',z'$
be the vertices corresponding to $y,z$ in $H(u_1)$, $x',z''$ be the
vertices corresponding to $x,z$ in $H(u_i)$ and $x'',y''$ be the
vertices corresponding to $x,y$ in $H(u_j)$. Since $\tau_3(H)=\ell$,
it follows that there exist $\ell$ internally disjoint pedant
Steiner trees connecting $\{x,y',z'\}$ in $H(u_1)$, say
$T_{1}',T_2',\cdots,T_{\ell}'$. Note that $\bigcup_{i=1}^k T_i$ is a
subgraph of $G$, $\bigcup_{j=1}^\ell T_j'$ is a subgraph of $H$, and
$(\bigcup_{i=1}^k T_i)\Box (\bigcup_{j=1}^\ell T_j')$ is a subgraph
of $G\Box H$. From Proposition \ref{pro2-5}, for any $T_i,T_j \
(1\leq i\neq j\leq k)$ and any $T_r',T_s' \ (1\leq r\neq s\leq
\ell)$, $(T_i\cup T_j)\Box (T_r\cup T_s)$ contains internally
disjoint pedant $S$-Steiner trees. Since $k\leq \ell$, there exist
$3\lfloor k/2 \rfloor$ internally disjoint pedant $S$-Steiner trees
in $(\bigcup_{i=1}^k T_i)\Box (\bigcup_{j=1}^\ell T_j')$, and hence
there are $3\lfloor k/2 \rfloor$ internally disjoint pedant
$S$-Steiner trees in $G\Box H$.

From the above argument, we conclude, for any $S\subseteq
V(G\Box H)$, that
$$
\tau_{G\Box H}(S)\geq \tau_{(\bigcup_{i=1}^k T_i)\Box (\bigcup_{j=1}^\ell T_j')}(S)\geq 3\lfloor \ell/2 \rfloor,
$$
which implies that $\tau_3(G\Box
H)\geq
3\lfloor \ell/2 \rfloor=3\lfloor \tau_3(H)/2 \rfloor$. The proof is complete.\qed

\section{Applications}

In this section, we demonstrate the usefulness of the proposed
constructions by applying them to some instances of Cartesian
product networks.

Given a vertex $x$ and a set $U$ of vertices, an \emph{$(x,U)$-fan}
is a set of paths from $x$ to $U$ such that any two of them share
only the vertex $x$. The size of an $(x,U)$-fan is the number of
internally disjoint paths from $x$ to $U$.

\begin{lem}{\upshape (Fan Lemma, \cite{West}, p-170)}\label{lem3-1}
A graph is $k$-connected if and only if it has at least $k+1$
vertices and, for every choice of $x$, $U$ with $|U|\geq k$, it has
an $(x,U)$-fan of size $k$.
\end{lem}

Hager also obtained the following result.

\begin{lem}{\upshape \cite{Hager}}\label{lem3-2}
Let $G$ be a graph, and let $k$ be an integer with $k\geq 2$. Then
$$
\tau_k(G)\geq \tau_{k+1}(G).
$$
\end{lem}

In \cite{Spacapan}, \u{S}pacapan obtained
the following result.

\begin{lem}{\upshape\cite{Spacapan}}\label{lem3-3}
Let $G$ and $H$ be two nontrivial graphs. Then
$$
\kappa(G\Box
H)=\min\{\kappa(G)|V(H)|,\kappa(H)|V(G)|,\delta(G)+\delta(H)\}.
$$
\end{lem}

\subsection{Two-dimensional grid graph and $n$-dimensional mesh}

A \emph{two-dimensional grid graph} is an $m\times n$ graph
$G_{n,m}$ that is the Cartesian product $P_n\Box P_m$ of path graphs
on $m$ and $n$ vertices. For more details on grid graph, we refer to
\cite{Calkin, Itai}.

\begin{pro}\label{pro3-4}
Let $n$ and $m$ be two integers with $n\geq 3, m\geq 3$. The network
$P_n\Box P_m$ has no pedant Steiner tree
connecting any three nodes. The number of internally disjoint pedant
Steiner trees is the maximum.
\end{pro}
\begin{pf}
From Theorem \ref{th1-4}, we have $\tau_3(P_n\Box P_m)\geq
\lfloor\frac{\tau_3(P_n)}{2}\rfloor+\lfloor\frac{\tau_3(P_m)}{2}\rfloor=0$.
Choose a vertex of degree $2$ in $P_n\Box P_m$, say $x$. Let $y,z$
be two neighbors of $x$. Then there is no internally disjoint pedant
Steiner tree connecting $\{x,y,z\}$. Therefore, $\tau_3(P_n\Box
P_m)=0$. \qed
\end{pf}\\

An \emph{$n$-dimensional mesh} is the Cartesian product of $n$
paths. By this definition, two-dimensional grid graph is a
$2$-dimensional mesh. An $n$-dimensional hypercube is a special case
of an $n$-dimensional mesh, in which the $n$ linear arrays are all
of size $2$; see \cite{Johnsson}.

\begin{cor}\label{pro3-5}
Let $k$ be a positive integer with $k\geq 3$. For $n$-dimensional
mesh $P_{m_1}\Box P_{m_2}$ $\Box \cdots \Box P_{m_n}$,
$$
\tau_k((P_{m_1}\Box P_{m_2}\Box \cdots \Box P_{m_n})=0.
$$
\end{cor}
\begin{pf}
$(1)$ From Proposition \ref{pro3-4}, $\tau_3((P_{m_1}\Box
P_{m_2}\Box \cdots \Box P_{m_n})=0$, and hence
$$
\tau_k((P_{m_1}\Box P_{m_2}\Box \cdots \Box P_{m_n})=0
$$
by Lemma \ref{lem3-2}.\qed
\end{pf}

\subsection{$n$-dimensional torus}

An \emph{$n$-dimensional torus} is the Cartesian product of $n$
cycles $C_{m_1},C_{m_2},\cdots,C_{m_n}$ of size at least three. The
cycles $C_{m_i}$ are not necessary to have the same size. Ku et al.
\cite{Ku} showed that there are $n$ edge-disjoint spanning trees in
an $n$-dimensional torus.

\begin{pro}\label{pro3-6}
$(1)$ For network $C_{m_1}\Box C_{m_2}$,
$$
\tau_k(C_{m_1}\Box C_{m_2})=
\begin{cases}
1,&\mbox k=3;\\
0,&\mbox k\geq 4,
\end{cases}
$$
where $m_i$ is the order of $C_{m_i}$ and $1\leq i\leq 2$.

$(2)$ Let $k$ be a positive integer with $k\geq 3$. For network
$C_{m_1}\Box C_{m_2}\Box \cdots \Box C_{m_n}$,
$$
0\leq \tau_k(C_{m_1}\Box C_{m_2}\Box \cdots \Box C_{m_n})\leq
2n-k+2,
$$
where $m_i$ is the order of $C_{m_i}$ and $1\leq i\leq n$.
\end{pro}
\begin{pf}
$(1)$ Set $G=C_{m_1}\Box C_{m_2}$. Since $\delta(G)=3$, it follows
from Lemma \ref{lem1-1} and Lemma \ref{lem3-2} that $\tau_3(G)\leq
1$ and $\tau_k(G)=0$ for $k\geq 4$. Since $\kappa(G)=3$, there
exists an $(x,S)$-fan for any $S\subseteq V(G)$ and $|S|=3$, where
$x\in V(G)\setminus S$. Then we have $\tau(S)\geq 1$ for any
$S\subseteq V(G)$ and $|S|=3$, and hence $\tau_3(G)=1$.

$(2)$ From Lemma \ref{lem3-3}, we have $\kappa(C_{m_1}\Box
C_{m_2}\Box \cdots \Box C_{m_n})=2n$, and hence
$$
0\leq \tau_k(C_{m_1}\Box C_{m_2}\Box \cdots \Box C_{m_n})\leq 2n-k+2
$$
by Lemma \ref{lem1-2}.\qed
\end{pf}

\subsection{$n$-dimensional generalized hypercube and $n$-dimensional hyper Petersen network}

Let $K_m$ be a clique of $m$ vertices, $m\geq 2$. An
\emph{$n$-dimensional generalized hypercube} \cite{DayA,
Fragopoulou} is the product of $m$ cliques. We have the following:

\begin{pro}\label{pro3-7}
Let $k$ be a positive integer with $k\geq 3$. For network
$K_{m_1}\Box K_{m_2}\Box \cdots \Box K_{m_n}$ where $m_i\geq k \
(1\leq i\leq n)$,
$$
\tau_k(K_{m_1}\Box K_{m_2}\Box \cdots \Box K_{m_n})\leq
\sum_{i=1}^{n}m_i-n-k+2.
$$
\end{pro}
\begin{pf}
From Lemma \ref{lem3-3}, we have $\kappa(K_{m_1}\Box K_{m_2}\Box
\cdots \Box K_{m_n})=\sum_{i=1}^{n}m_i-n$, and hence
$$
0\leq \tau_k(C_{m_1}\Box C_{m_2}\Box \cdots \Box C_{m_n})\leq
\sum_{i=1}^{n}m_i-n-k+2
$$
by Lemma \ref{lem1-2}.\qed
\end{pf}\\

An \emph{$n$-dimensional hyper Petersen network} $HP_n$ is the
product of the well-known Petersen graph and $Q_{n-3}$ \cite{Das},
where $n\geq 3$ and $Q_{n-3}$ denotes an $(n-3)$-dimensional
hypercube. Note that $HP_3$ is just the Petersen graph.

\begin{pro}\label{pro3-8}
$(a)$ The network $HP_3$ has one pendant Steiner tree
connecting any three nodes.

$(b)$ The network $HP_4$ has two internally disjoint pedant Steiner trees
connecting any three nodes. The number of internally disjoint pendant Steiner trees
is the maximum.
\end{pro}
\begin{pf}
$(a)$ Note that $HP_3$ is just the Petersen graph. Set $G=HP_3$.
Since $\delta(G)=3$, it follows that $\tau_3(G)\leq 1$ by Lemma
\ref{lem1-1}. From Lemma \ref{lem3-2}, there exists an
\emph{$(x,S)$-fan} for any $S\subseteq V(G)$ and $|S|=3$, where
$x\in V(G)\setminus S$. Thus $\tau(S)\geq 1$, and hence
$\tau_3(G)=1$, that is, $HP_3$ has one pedant Steiner tree
connecting any three nodes.

$(b)$ Since $\delta(G)=4$, it follows from Lemma \ref{lem1-1} that
$\tau_3(HP_4)\leq 2$. One can check that for any $S\subseteq V(G)$
and $|S|=3$, $\tau(S)\geq 2$. So $\tau_3(G)=2$. \qed
\end{pf}

\end{document}